\documentclass[11pt,epsf]{article}
\usepackage{graphicx}
\usepackage{theorem, amsmath, amssymb}
\newtheorem{Theorem}{Theorem}[section]
\newtheorem{Definition}[Theorem]{Definition}
\newtheorem{Proposition}[Theorem]{Proposition}

\newtheorem{Lemma}[Theorem]{Lemma}
\newtheorem{Corollary}[Theorem]{Corollary}

\theoremstyle{remark}
\newtheorem{Remark}[Theorem]{Remark}

  {\begin{equation}\left\{\begin{aligned}}%
  {\end{aligned}\right.\end{equation}\ignorespacesafterend}%

\newenvironment{SDE*}%
  {\begin{equation*}\left\{\begin{aligned}}%
  {\end{aligned}\right.\end{equation*}\ignorespacesafterend}%

\usepackage{color}

\begin{document}
\title{On  martingale problems 
with continuous-time mixing 
and values of zero-sum games 
 without Isaacs condition
}

\author{Mihai S\^{\i}rbu \footnote{University of Texas at Austin,
    Department of Mathematics, 1 University Station C1200, Austin, TX,
    78712.  E-mail address: sirbu@math.utexas.edu. The research of
    this author was supported in part by the National Science
    Foundation under Grant    DMS 1211988. Any opinions, findings, and conclusions or recommendations expressed in this material are those of the authors and do not necessarily reflect the views of the National Science Foundation.}}
\maketitle
\begin{abstract} We consider  a zero-sum stochastic differential game over  elementary mixed feed-back strategies. These are strategies based only on the knowledge of the past state,  randomized continuously in time from a sampling distribution which is kept constant in between some stopping rules. Once both players choose such strategies, the state equation admits a unique solution in the sense of the martingale problem of Stroock and Varadhan. We show that  the game  defined over  martingale solutions has a value, which is the unique continuous viscosity solution of  the randomized Isaacs equation.
\end{abstract}


\section{Introduction}Continuous-time zero-sum  stochastic differential games have been studied extensively since \cite{fs}. If the Isaacs condition holds, one expects the game to have a value, which is the solution of the (unified) Isaacs equation. Many interesting results in this direction have been obtained. However, 
to the best of our knowledge, the only work dealing with the existence of a value for general continuous-time zero-sum stochastic differential games \emph{without Isaacs conditions} is the very recent work \cite{buc-li-q}. 
The authors consider a game where players see each other's actions with a delay relative to a fixed time grid, and both play mixed delayed strategies. The cost functional is defined through  a BSDE, following the earlier work \cite{hl}.  Letting the time-grid become finer, the lower and the upper values of the approximate games converge to the \emph{same} value, which  is the unique viscosity solution of the randomized Isaacs equation. Since the time-grid is first fixed, this is resembling of  \cite{MR945913}, but with randomization. The  very  definition of the game and mixed strategies  is highly non-trivial and the analysis is technical.  Deterministic games over (discretized)  mixed feed-back strategies (called positional/historical)  without  the presence of Isaacs conditions  were studied in \cite{krasovskii-subbotin-88} using different methods.

 \noindent  To summarize, the program we propose here for continuous-time stochastic games without Isaacs condition  is to:
\begin{enumerate}
\item model the zero-sum game symmetrically, over generalized (in the sense of \cite{young}) feedback strategies, that we interpret as mixed continuously in time (see Definition \ref{def:gen-feedback} and the subsequent comments),

\item  posit the state equation in the sense of the martingale problem (similar to \cite{ek-compactness} in case of one player with open-loop controls)  but only on the canonical space of $X$ (as in  Stroock and Varadhan  \cite{MR2190038}, followed by  \cite{pz-game}),
\item identify a simple class (called elementary) of generalized feed-back strategies over which the state equation is well-posed. This is actually an important task if one chooses to work with feedback strategies  (see the comments in  \cite{elliott-feedback}) and whose resolution is non-trivial in our general model. In the case where only the drift is controlled by the two players, the well posed-ness over feedback strategies can be studied using Girsanov's theorem (this is the case, for \emph{general} pure feed-back strategies, in \cite{elliott-stochastic-games}) but the same cannot be applied here.
\item show that the game does have a value over the class of strategies identified above, and the value is the unique continuous viscosity solutions of the Isaacs equation. This is done by reducing the problem to previous work \cite{sirbu} over pure feed-back strategies.  Obviously, there exist approximate saddle points. In addition, we show that a version of the Dynamic Programming Principle holds.
\end{enumerate}
In order to accomplish Items 3 and 4 (together) we make use of the deep theory of Yamada-Watanabe relating strong existence and path-wise uniqueness to weak existence and uniqueness for functional  SDE's.
We emphasize right away that, having a well-posed state equation in Item 3 above requires a class of strategies that does not contain (even formally) a saddle point, unlike the case when only the drift is controlled in \cite{elliott-stochastic-games}.

We attempt to better relate our program to existing work below.
Following earlier  stochastic control literature (see the notion of natural strategies in  \cite{MR2723141})  and  the work on deterministic games in  \cite{elliott-feedback}, \cite{krasovskii-subbotin-88} or  on stochastic games in \cite{elliott-stochastic-games}, 
we   consider here the game over feed-back strategies, as in the more recent papers \cite{sirbu} or \cite{pz-game}. However, our strategies are randomized (sampled) \emph{continuously} in time. Actually (in the case of open-loop controls customary in one-player/control problems), continuous-time randomization  amounts to choosing relaxed/generalized controls in the sense of L.C.Young \cite{young} (see the comments after Definition \ref{def:gen-feedback}). The concept of  open-loop generalized/relaxed controls  has been first used in  the context of games  in \cite{smoljakov} (for a particular example of deterministic games) and later in \cite{elliott-kalton-markus} to obtain existence of a saddle point for a linear deterministic game without Isaacs condition.  We have not seen stochastic games posed over generalized feedback strategies, the way we use here.

In our choice of feed-back strategies, in between stopping times based only on the knowledge of the past of the state, the player samples continuously and independently in time from the same distribution of possible actions,  and the sampling distribution is chosen at the earlier time depending on the past of the state.
One does not  expect the state equation to have strong solutions in such a formulation  over mixed feed-back strategies, and the continuous-time mixing is not even explicitly modeled.
The only information needed from the state equation is the law of the state process, so 
 we study the  existence and uniqueness of solutions to the state equation  in the sense of the martingale problem of Stroock and Varadhan  \cite{MR2190038}.  Following \cite{ek-compactness},  in   martingale  formulation, the continuous-time randomization can be averaged out  similarly to the way generalized controls in \cite{young} average out locally in time leading to generalized curve solutions of the state equation.
 Due to the fact that we work with feed-back strategies, our formulation of the martingale problem is posed over the space of paths of the controlled process, like in \cite{pz-game} (and following  \cite{MR2190038}), as opposed to the more classic literature on stochastic control where the the canonical space has to accommodate the open loop controls (relaxed or not) like in \cite{ek-compactness}.


Compared to \cite{buc-li-q}, we work with the original game, without any discrete-time approximation/exogenous restriction, but we allow for continuous-time mixing/generalized strategies and (for that reason) consider martingale solutions of the state equations. Aside from considering the game over feed-back strategies, we may say that we substitute the limit over time-grids in \cite{buc-li-q} by the very definition of solutions to the state equations. We believe this provides some more understanding to continuous time games with randomization. In particular,  the ``weaker'' player (the one in the exterior of inf/sup or sup/inf) needs to randomize continuously in time  i.e. use feed-back strategies   which are relaxed/generalized in the sense of \cite{young}.  Technically, we also allow for local Lipschitz conditions rather than global Lipschitz  in the description of the game.

\section{The set-up and main results}
In what follows, we borrow heavily the notation and hypotheses from \cite{sirbu}.  { In order to make the presentation self-contained, a few arguments in \cite{sirbu} are repeated, but the overlap is kept to a minimum.
We consider a stochastic differential game with two players.
First player's actions belong to a compact
$ U\subset \mathbb{R}^k$, while second player's action belong also to a compact set  $V\subset \mathbb{R}^l$. The state lives in  $\mathbb{R}^d$.
Let $b:[0,T] \times \mathbb{R}^d \times U \times V \to \mathbb{R}^d$ and $\sigma:[0,T]\times\mathbb{R}^d\times U \times V \to \mathbb{M}_{d,d'}$ be two continuous functions. We consider the diffusion 
\begin{equation}\label{eq:SDE}\left\{
\begin{array}{ll}
 dX_t=b(t,X_t,u_t, v_t)dt+\sigma (t, X_t,u_t, v_t)dW_t,   \\ 
 X_s=x \in \mathbb{R}^d,
 \end{array} \right.
 \end{equation}
 starting at initial time $s$ at position $x$, and 
which is controlled by both players. Here, $W$ is a $d'$-dimensional Brownian motion on a \emph{fixed} probability space
 $(\Omega, \mathcal{F}, \mathbb{P})$. Consider the natural augmented filtration
 $\mathbb{F}^s=(\mathcal{F}^s_t)_{s\leq t\leq T},$ generated by the Brownian increments starting at $s$, i.e., 
 $$\mathcal{F}^s_t=\sigma (W_u-W_s, s\leq u\leq t)\vee \mathcal{N}(\mathbb{P}, \mathcal{F})\ \  \textrm{for}\ \ s\leq t\leq T.$$
 Now, given a   bounded and continuous function $g:\mathbb{R}^d\rightarrow \mathbb{R}$, the second player pays, in this game, to the first player, the amount 
 $\mathbb{E}[g(X^{s,x;u,v}_T)]$ leading to the  zero-sum game  
$$\sup _{u }\inf _v\mathbb{E}[g(X^{s,x;u,v}_T)],\ \ \ \inf _v\sup _u  \mathbb{E}[g(X^{s,x;u,v}_T)].\ \ \ \ \  \ 
$$
We denote by $\mathcal{P}(U)$ and $\mathcal{P}(V)$ all  probability measures on $(U, \mathcal{B}(U))$, respectively $(V,\mathcal{B}(V)$. Considered with the L\'evy-Prokhorov metric (the metric of weak convergence of probability measures, since $U,V$ are separable), these are compact metric spaces also.
We make the notation
$$L(t,x,u,v, p,M)\triangleq \left[b(t,x,u,v)\cdot p+\frac{1}{2}Tr \left (\sigma(t,x,u,v)\sigma(t,x,u,v)^T M\right)\right].$$
 It is well known (see, for example, relation (1.11) in \cite{buc-li-q}) that, since the spaces $\mathcal{P}(U)$ and $\mathcal{P}(V)$ are convex and compact and 
$(\mu,\nu)\rightarrow \int _{U\times V}
L(t,x, u,v, p, M)\mu (du) \nu (dv)$ is a bi-linear form, under minimal  continuity hypotheses, one can define 
the randomized  Hamiltonian as 
\begin{eqnarray*}
H(t,x,p,M)\triangleq \sup_{\mu \in   \mathcal{P}(U)} \ \ \ \inf _{\nu \in \mathcal{P}(V) } \
\int _{U\times V}
L(t,x, u,v, p, M)\mu (du) \nu (dv)\\
= \ \ \ \inf _{\nu \in \mathcal{P}(V) } \ \  \sup_{\mu \in   \mathcal{P}(U)}\int _{U\times V}
L(t,x, u,v, p, M)\mu (du) \nu (dv).
\end{eqnarray*}
In other words, the equality above holds even if the Isaacs condition does not hold, i.e  even if \\
$$H^-(t,x,p,M)\triangleq \sup_{u \in   U} \ \ \ \inf _{v\in V } 
L(t,x, u,v, p, M)<
\inf _{v \in  V } \ \  \sup_{u \in U }
L(t,x, u,v, p, M) \triangleq H^+(t,x,p,M).
$$
We associate the randomized  Isaacs equation to the game
\begin{equation}\label{eq:Isaacs}
\left \{
\begin{array}{ll}
-v_t-H (t,x,v_x,v_{xx})=0\ \ \textrm{on}\ [0,T)\times \mathbb{R}^d,\\
v(T,\cdot)=g(\cdot),\ \ \textrm{on}\ \mathbb{R}^d.
\end{array}
\right.
\end{equation}

\noindent {\bf Standing Assumptions:}
the coefficients $b, \sigma$ of the  stochastic system
\begin{enumerate}
\item  are jointly continuous in $(t,x,u,v) \in [0,T]\times \mathbb{R}^d\times U\times V$
\item  satisfy a local but uniform Lipschitz condition in $x$, i.e.
$$
{\bf(L)} \ |b(t,x,u,v)-b(t,y,u,v)|+|\sigma (t,x,u,v)-\sigma (t,y,u,v)|\leq L(K) |x-y|\ \forall\  |x|,|y|\leq K,$$
$ \forall \  t\in [0,T],  u\in U, v\in V$  for some $L(K)<\infty$
\item satisfy a global linear growth condition, i.e. there exists $C<\infty$ such that 
$${\bf (LG)}  \ |b(t,x,u,v)|+|\sigma (t,x,u,v)|\leq C(1+|x|), \forall\  x, y\in \mathbb{R}^d, t\in [0,T],  u\in U, v\in V$$ 
\item either $d=1$ (one-dimensional state) or $\sigma$ is smooth in $x$, i.e. for each fixed $t,u,v$ there exists 
$\frac{\partial ^2 \sigma}{\partial x_i x_j} (t,x,u,v)$  continuous in $x\in \mathbb{R}^d$, for all $i,j=1, d$ and uniformly locally bounded in $x$, i.e.
$$\left |\frac{\partial ^2 \sigma}{\partial x_i x_j}(t,x,u,v)\right|\leq C(K),\ \ \forall |x|\leq K, t\in [0,T],  u\in U, v\in V.$$
\end{enumerate}
Using only items 2 and 3 from the Standing Assumption, we see that  if both players choose to hold a constant strategy, the controlled state has a unique solution. This holds even if the game is started at a random (but stopping) time sequel to the starting time $s$, and the constant strategies depend on whatever happened before the starting time.  More precisely, we have the following proposition:
\begin{Proposition} \label{prop:state-eq}
\begin{enumerate}
\item Let $s\leq \tau' \leq T$ be a stopping time of the filtration $\mathbb{F}^s=(\mathcal{F}^s_t)_{s\leq t\leq T}$. Let 
$\xi \in \mathbb{R}^d, a\in U, b \in V, $ be random variables measurable with respect to $\mathcal{F} ^s_{\tau '}.$ Then, the system starting at time $\tau '$ with initial condition $\xi$, where both players choose constant strategies $a,b$ in between $\tau '$ and $T$ has a unique strong solution. In other words, the SDE
\begin{equation}\label{state-simple}
\left \{
\begin{array}{ll}
dX_t=b(t,X_t, a, b)dt+\sigma (t, X_t,a, b)dW_t,\ \tau ' \leq t\leq T\\
X_{\tau '}=\xi
\end{array}
\right .
\end{equation}
has a unique  strong solution $(X_t)_{\tau ' \leq t\leq T}$. 
\item In addition, we have path-wise uniqueness. More precisely, consider  any probability space accommodating the same Brownian motion on the time interval $[s,T]$, with respect to two different filtrations  $\mathbb{F}_i^s=(\mathcal{F}^i_t)_{s\leq t\leq T}$ satisfying the usual conditions. If one chooses an initial stopping time $\tau'$ with respect to both filtrations, and $a,b\in \mathcal{F}^i_{\tau '}$ for $i=1,2$, any two solutions $(X^i_t)_{\tau ' \leq t\leq T}$ of \eqref{state-simple}corresponding to the two filtrations for $i=1,2$ satisfy
$$\mathbb{P}(X^1_t=X^2_t, \forall \ \tau '\leq t\leq T)=1.$$
\end{enumerate}

\end{Proposition}
Proof: the proof is standard (and even in \cite{sirbu} most details were omitted). $\diamond$
 
For the remainder of the paper, fixed a starting time $s$, we denote by $C([s,T])\triangleq C([s,T],\mathbb{R}^d)$ and endow this path space with the natural (and raw) filtration 
 $\mathbb{B}^s=(\mathcal{B}^s_t)_{s\leq t\leq T}$ defined by
 $$\mathcal{B}^s_t\triangleq \sigma (y(u),s\leq u\leq t), \ \ s\leq t\leq T.$$
 Elements of the path space $C([s,T])$ will be denoted by $y (\cdot)$ or $y$. Stopping times on the space $C([s,T])$ with respect with the filtration $\mathbb{B}^s$, i.e.  mappings $\tau :C([s,T])\rightarrow [s,T]$ satisfying
 $  \{\tau \leq t\}\in \mathcal{B}^s_t \ \forall \ s\leq t\leq T$
are called  stopping  rules, following \cite{ks}. We denote   by $ \mathbb{B}^s$ the class of such stopping  rules starting at $s$.

 \begin{Definition}[Pure   Feed-Back Strategies]  \label{def:gen-feedback-pure}Fix $0\leq s\leq T$.
 A pure feed-back  strategy $u$ starting at $s$,  for the first player,  is  a mapping
  $u: (s,T]\times C([s,T])\rightarrow U$  which is predictable with respect to the raw filtration  $\mathbb{B}^s$.  
    A  pure strategy $v$ for the second player is defined in an identical way, but takes values in $V$.
 We denote by $\overline{\mathcal{U}} ^p(s)$ and $\overline{\mathcal{V}}^p(s)$ the collections of all possible pure strategies for the first, and the second player, given the initial deterministic time $s$.
\end{Definition}

Next definition is borrowed from \cite{sirbu}.
 \begin{Definition}[Elementary Pure   Feed-Back Strategies] \label{def:s} Fix $0\leq s\leq T$.
 An   elementary pure feed-back  strategy $u$ starting at $s$,  for the first player,  is defined by

\begin{itemize}
\item a finite non-decreasing sequence of stopping  rules, i.e. $\tau _k \in \mathbb{B}^s$ for $k=1, \dots, n$ and 
$$s=\tau _0\leq \dots \tau _k\leq \dots \leq \tau _n=T $$
\item for each $k=  1\dots n$, a constant value of the strategy $\xi_k$ in between the times $\tau _{k-1}$ and $\tau _k$, which is decided based only on the knowledge of the past state up to  $\tau _{k-1}$, i.e.
$\xi_k:C([s,T])\rightarrow U$ such that
$\xi _k\in \mathcal{B}^s_{\tau _{k-1}}$.
\end{itemize}
 The strategy is to hold $\xi_k$ in between $(\tau _{k-1}, \tau _{k}]$, i.e.
 $u: (s,T]\times C([s,T])\rightarrow U$ is defined by
 $$u(t, y(\cdot))\triangleq \sum _{k=1}^n \xi _{k}(y(\cdot))1_{\{ \tau _{k-1}(y(\cdot))<t\leq  \tau _{k}(y(\cdot))\}}.$$
  An elementary pure strategy $v$ for the second player is defined in an identical way, but takes values in $V$.
 We denote by $\mathcal{U} ^p(s)$ and $\mathcal{V}^p(s)$ the collections of all possible pure elementary strategies for the first, and the second player, given the initial deterministic time $s$.
\end{Definition} 
Simply iterating Proposition \ref{prop:state-eq} we obtain (already seen in 
\cite{sirbu}, actually) that 
the  elementary pure  strategies we just defined produce strong solutions of the state equation without any  Lipschitz assumption in $u,v$.   The same is not true for  the general feed-back strategies in Definition \ref{def:gen-feedback-pure}.
\begin{Proposition}\label{prop:state-eq-simple} Fix $s,x$ and 
  let  players one and two choose strategies $u\in \mathcal{U}^p(s)$ and $ v \in\mathcal{V}^p(s)$. Then, there exists a unique strong (and square integrable) solution
$(X^{s,x;u,v}_t)_{s \leq t\leq T}$, $X^{s,x;u,v}_t\in \mathcal{F}^s_t$ of the state equation
\begin{equation}\label{state-full}
\left \{
\begin{array}{ll}
dX_t=b(t,X_t,  u (X_{\cdot}), v(X_{\cdot}))dt+\sigma (t, X_t,u (X_{\cdot}), v(X_{\cdot}))dW_t,\ s \leq t\leq T\\
X_{s}=x\in \mathbb{R}^d.
\end{array}
\right .
\end{equation} 
In addition, the (functional) SDE \eqref{state-full} satisfies the path-wise uniqueness property.
\end{Proposition}
The restriction  to the class of elementary feed-back strategies is motivated both by modeling   considerations (actions can be changed discretely in time) and by technical reasons  (the state equation is well-posed in a strong sense).
It is now time to introduce the strategies which will be used for the game in the absence of Isaacs condition.
 \begin{Definition}[Continuous-Time Mixed  Feed-Back Strategies] \label{def:gen-feedback}
Fix $s$. 
A  continuous-time randomized feed-back strategy $\mu$  for the first player, 
 starting at $s$,   is a mapping
 $\mu : (s,T]\times C([s,T])\rightarrow \mathcal{P}(U)$  which is predictable with respect to the raw filtration $\mathbb{B}^s$.
  A   continuous-time mixed strategy $\nu $ for the second player is defined in an identical way, but takes values in $\mathcal{P}(V)$.
 We denote by $\overline{\mathcal{U}}^m (s)$ and $\overline{\mathcal{V}}^m(s)$ the collections of all possible  elementary  continuous time mixed strategies for the first, and the second player, given the initial deterministic time $s$.\end{Definition}
We have defined our  mixed feed-back strategies  simply as 
$\mu (\nu): (s,T]\times C([s,T])\rightarrow \mathcal{P}(U)(\mathcal{P}(V)),$ 
 but it is understood that, given two strategies $\mu\in \overline{\mathcal{U}}^m(s)$ and $\nu\in \overline{\mathcal{V}}^m(s)$, the two players will sample, continuously in time, independently from the past and from each other, from the sampling distributions
 $\mu (t,y(\cdot);du)$ and, respectively,  $\nu (t,y(\cdot);du)$.
As mentioned in the introduction, in the case of open-loop controls (not present here),  what we call continuous-time mixing is similar to  choosing relaxed controls in the sense of \cite{young}. It represents a very particular and extreme form of mixing.
Following the use of relaxed open-loop controls in the seminal work of \cite{young} (deterministic control) and later \cite{ek-compactness} (stochastic control), 
we do not attempt to define here \emph{explicitly} the continuous-time mixing mechanism which is used, to sample at each $t$, independently from any other time and from the other player, from the corresponding sampling distribution,  but  describe it  through a martingale problem  involving averaging over $u$ and $v$ of the generator, locally in time.

In  the case of pure strategies,  the existence of strong solutions for the state equation is, in our opinion, important: if both players have decided precisely what to do contingent on the past state, then the noise should have a one-to-one response through the state equation.   In some sense,  a similar argument is made for deterministic games in \cite{elliott-feedback}.
 In the case of mixing, one needs additional randomness anyway to define a solution of the state equation, and the additional randomness is chosen by the players. 
 We therefore  study the law of the state process $X^{s,x,\mu,\nu}$ via the well posedness  of the martingale problem of Stroock and Varadhan. Since we are using feed-back strategies, the martingale problem can be posed directly on the canonical space of the state-path space, as opposed to needing a larger space to accommodate open-loop relaxed controls in the one-player games in \cite{ek-compactness}.
A similar  martingale formulation is considered  in \cite{pz-game}  but since pure strategies are considered, no averaging over mixing is needed in the definition.

 We would like to point out that the definition below only makes sense for feed-back strategies, and one would have a harder time posing such a problem in an Elliott-Kalton formulation of the game. Also, the very definition below replaces, in some conceptual sense, the limiting arguments in \cite{buc-li-q}, as mentioned in the Introduction.
\begin{Definition}[Martingale Solutions over Mixed Strategies]\label{mart} Fix $s,x$ and $\mu \in \overline{\mathcal{U}}^m(s)$, $\nu \in \overline{\mathcal{V}}^m(s)$. A probability measure $\mathbb{Q}$ on $(C[s,T], \mathcal{B}(C[s,T]))$ solves the martingale problem for $s,x,\mu,\nu$ if, for each $f\in C^2(\mathbb{R}^d)$ we have that
the process $(M^f_t)_{s\leq t\leq T}$ defined by 
$$M^f_t(y)\triangleq f(y(t))-\int _s^t \left ( \int _{U\times V}L(r, y(r), u, v, f_x (r,y(r)), f_{xx}(r, y(r))) \mu (r, y(\cdot); du)\nu (r, y(\cdot); dv) 
\right ) dr $$
is a continuous local martingale on $(C[s,T], \mathbb{B}^s, \mathbb{Q})$ and $\mathbb{Q}(y(s)=x)=1$.
\end{Definition}
The heuristics behind Definition \ref{mart} is apparent,  and  shows that our interpretation of continuous-time mixing  is similar to the idea of relaxed controls in \cite{young}, later  used in  martingale formulation of open-loop control problems over relaxed controls in in \cite{ek-compactness}. More precisely, if the two players choose mixed strategies $\mu$ and $\nu$, then at time $t$, they both sample independently from each other and from whatever happened up to time $t$ from $\mu (t, y(\cdot); du)$  respectively  $\nu (t, y(\cdot); dv)$ to choose their actions for the infinitesimal future. This means that, abusing notation (among others since the ``past'' has not been precisely defined) we formally have, 
$$\mathbb{E}_t[X_{t+dt}-X_t]= \left (\int _{U\times V}b(t, y(t), u, v) \mu (t, y(\cdot); du)\nu (t, y(\cdot); dv)\right) dt $$
and
$$Var_t(X_{t+dt}-X_t)= \left (\int _{U\times V}\sigma (t, y(t), u, v)  \sigma ^T(t, y(t), u, v) \mu (t, y(\cdot); du)\nu (t, y(\cdot); dv)\right )dt. $$
Therefore, the probability measure $\mathbb{Q}$ in the above definition has the obvious meaning of the law of the state process.
We point out that posing the (local) martingale problem over the raw filtration $\mathbb{B}^s$ on $C[s,T]$ is equivalent to posing the martingale problem over $\overline{\mathbb{B}^s} ^{\mathbb{Q}}$, the smallest filtration satisfying the usual conditions under $\mathbb{Q}$ and containing $\mathbb{B}^s$. 
We identify $u\in U\equiv \delta _u\in \mathcal{P}(U)$, resulting in the continuous embedding $U\subset \mathcal{P}(U).$ Therefore, we also have the  embeddings
$$\mathcal{U}^p(s)\subset \overline{\mathcal{U}}^p(s)\subset  \mathcal{U}^m(s), \mathcal{V}^p(s)\subset \overline{\mathcal{V}}^p(s)\subset \mathcal{V}^m(s).$$
{If the game is deterministic ($\sigma=0$), using strategies which are randomized continuously in time would formally result in a deterministic state equation. In light of the Strong Law of Large Numbers (continuous version), this is  not surprising and was well understood since \cite{krasovskii-icm-78} in the context of games of feed-back strategies, and well before that in the context of open-loop relaxed controls in \cite{young}.

\begin{Remark} If $u\in \mathcal{U}^p(s), v\in \mathcal{V}^p(s)$ are elementary pure strategies, then, the law of the (path-wise) unique solution $X^{s,x,u,v}$ i.e.
$$\mathbb{Q}\triangleq \mathbb{P}\circ \left (X^{s,x,u,v}_{\cdot}\right)^{-1}$$
is the unique solution of the martingale problem. Uniqueness comes from two observations: the martingale problem is equivalent to the existence of a weak solution ( \cite{MR2190038} or \cite{KS88}, Corollary 4.8 page 317)) and path-wise uniqueness implies uniqueness in Law, from the celebrated result of Yamada-Watanabe (see \cite{KS88}, Proposition 3.20 on page 309). This program works well even for functional equations that we consider here, i.e. allows for dependence on the whole past of the path.
\end{Remark}
 We can define the  well-posedness of the state equation over the  class of  general mixed feed-back strategies in Definition \ref{def:gen-feedback}, but  we cannot expect such well-posedness to hold.  This is  one reason for which we do not study the existence of a saddle point at this level of generality. If a saddle point  exists formally (for example, as a limit in some sense of approximate saddle points), it could be a pair of  general mixed feed-back strategies  for which  the  state equation is not  well posed. Modeling a game over a class of mixed  feedback strategies which is large enough to contain saddle points, but small enough to make the state equation well posed is  most likely not possible in the general case. As pointed out in \cite{elliott-feedback} or \cite{krasovskii-subbotin-88} (even for deterministic games over feed-back strategies), it is important to model the game over a  (small enough) class of feed-back strategies that leads to a well-posed state equation. This is what we do here. 
 
Fortunately, the state equation is well posed in the martingale sense, over the  restricted class of elementary  mixed strategies,  defined below,  and we  only define the game  over such strategies.

\begin{Definition}[Elementary Continuous-Time Mixed  Feed-Back Strategies]Fix $s$. 
An elementary, but continuous-time randomized feed-back strategy $u$  for the first player, 
 starting at $s$,   is defined by  
\begin{itemize}
\item (again) a finite non-decreasing sequence of stopping  rules, i.e. $\tau _k \in \mathbb{B}^s$ for $k=1, \dots, n$ and 
$$s=\tau _0\leq \dots \tau _k\leq \dots \leq \tau _n=T $$
\item for each $k=  1\dots n$, a constant sampling distribution $\xi_k \in \mathcal{P}(U)$ in between the times $\tau _{k-1}$ and $\tau _k$, which is decided based only on the knowledge of the past state up to  $\tau _{k-1}$, i.e. a measurable 
$$\xi_k:(C([s,T]),  \mathcal{B}^s_{\tau _{k-1}})\rightarrow (\mathcal{P}(U),\mathcal{B}(\mathcal{P}(U))).$$ 
\end{itemize}
 The strategy is to sample continuously and independently in time from the sampling distribution  $\xi_k (du) $ in between $(\tau _{k-1}, \tau _{k}]$. More precisely, with 
 $\mu : (s,T]\times C([s,T])\rightarrow \mathcal{P}(U)$  defined by
 $$\mu (t, y(\cdot);du)\triangleq \sum _{k=1}^n \xi _{k}(y(\cdot))1_{\{ \tau _{k-1}(y(\cdot))<t\leq  \tau _{k}(y(\cdot))\}} \in \mathcal{P}(U), $$
the player samples, at each time t, from $\mu (t,y(\cdot);du)$ independently from any other time.
  An elementary  continuous-time mixed strategy $\nu $ for the second player is defined in an identical way, but takes values in $\mathcal{P}(V)$.
 We denote by $\mathcal{U}^m (s)$ and $\mathcal{V}^m(s)$ the collections of all possible  elementary  continuous time mixed strategies for the first, and the second player, given the initial deterministic time $s$.

\end{Definition}
Obviously, 
$\mathcal{U}^p(s)\subset  \mathcal{U}^m(s), \mathcal{V}^p(s)\subset \mathcal{V}^m(s).$
\begin{Proposition}\label{prop:well-posed} Let $\mu \in \mathcal{U}^m(s),\nu \in  \mathcal{V}^m(s) $ be two elementary continuous-time mixed strategies. There exists a unique solution $Q^{s,x,\mu,\nu }$ of the martingale problem in Definition \ref{mart}.
\end{Proposition}
 The restriction of mixed strategies $\overline{\mathcal{U}}^m/\overline{\mathcal{V}}^m$ to the elementary classes $\mathcal{U}^m/\mathcal{V}^m$, unlike  the case of pure strategies (in Proposition \eqref{prop:state-eq-simple} and in \cite{sirbu}),  has only a mathematical motivation,  namely to obtain a well-posed state equation. From the point of view of modeling, elementary strategies do not  have the meaning of  changing actions only at discrete times, as the randomization is still performed at any time, continuously.
The proof of the Proposition above is postponed until Section \ref{proofs}. Having this result,
we  can now formulate the game rigorously over elementary continuous time mixed  strategies. For  fixed $0\leq s \leq t \leq T$ and $x\in \mathbb{R}^d$ (deterministic) and $\mu \in \mathcal{U}^m(s)$ , $\nu \in \mathcal{V}^m(s)$, we define the expected amount that player two pays to player one
$$J(s,x,\mu ,\nu )\triangleq \mathbb{E}^{\mathbb{Q}^{s,x,\mu ,\nu }}[g(y(T))].$$ If $u\in \mathcal{U}^p(s)$ and $v\in \mathcal{U}^p(s)$ are elementary pure strategies, then, obviously,
$J(s,x,u,v)=\mathbb{E}[g(X^{s,x,u,v}_T)],$ { where $X^{s,x,u,v}$ is the (path-wise) unique strong solution of the state equation. For elementary mixed strategies $\mu\in \mathcal{U}^m, \nu\in \mathcal{V}^m$, the quantity $J(s,x,\mu,\nu)$ still has the meaning of $\mathbb{E}[g(X^{s,x,\mu,\nu}_T)]$. As previously mentioned, we do not address the existence and uniqueness of $X^{s,x,\mu,\nu}$ directly, but look only for what its law $\mathbb{Q}^{s,x,\mu,\nu}$ should be. In this respect, our set-up of the game resembles the framework of \cite{ms96} (going all the way back to  Blackwell), where the responses to players' actions are (local in time) distributions of the state, and not specific realizations of the state process.

We now  define the lower and the upper value of the game
$$
V^-(s,x)\triangleq \sup _{\mu  \in \mathcal{U}^m(s)} \inf _{\nu \in \mathcal{V}^m(s)} J(s,x,\mu ,\nu )\leq 
\inf _{\nu \in \mathcal{V}^m(s)}  \sup _{\mu  \in \mathcal{U}^m(s)}J(s,x,\mu ,\nu )\triangleq V^+(s,x).$$
The main result is
\begin{Theorem}
\label{main}
The game has a value, i.e. $V^-=V^+$ and the function 
$V\triangleq V^-=V^+$ is the unique bounded and continuous viscosity solution of the randomized Isaacs equation. In addition, we have the DPP
$$V(s,x)= \sup _{\mu \in \mathcal{U}^m(s)} \inf _{\nu \in \mathcal{V}^m(s)}  \mathbb{E}^{\mathbb{Q}^{s,x,\mu,\nu }}[V (\rho(y), y(\rho (y)))]=\inf _{\nu \in \mathcal{V}^m(s)}  \sup _{\mu \in \mathcal{U}^m(s)}  \mathbb{E}^{\mathbb{Q}^{s,x,\mu ,\nu }}[V (\rho(y), y(\rho (y)))], \forall \rho \in \mathbb{B}^s.$$
\end{Theorem}
As a simple corollary, we obtain
\begin{Corollary} For each $\varepsilon>0$ there exists an $\varepsilon$-saddle point, which means
$(\mu (\varepsilon),\nu (\varepsilon))\in \mathcal{U}^m\times \mathcal{V}^m$
such that 
$$J(s,x,\mu  ,\nu (\varepsilon)) -\varepsilon \leq J(s,x,\mu( \varepsilon),\nu (\varepsilon)) \leq J(s,x,\mu (\varepsilon) ,\nu) +\varepsilon\  (\forall) \ \ (\mu, \nu) \in \mathcal{U}^m\times \mathcal{V}^m.$$\end{Corollary}
\section{Proofs}\label{proofs}
The proofs  are short and are based on the simple idea to identify an auxiliary game, satisfying the Isaacs equation, whose solution of the state equation over pure strategies solves the martingale problem, and then appeal to \cite{sirbu}.
 This way we avoid the need to use some representation with respect to martingale measures in the spirit of \cite{ek-m-mart-meas}. More precisely, we define
$\tilde{b}:[0,T]\times \mathbb{R}^d\times \mathcal{P}(U)\times \mathcal{P}(V)\rightarrow \mathbb{R}^d$
by 
$$\tilde{b} (t,x, \mu, \nu)\triangleq \int _{U\times V}b(t,x,u,v)\mu (du) \nu (dv),$$ and 
$\tilde{\sigma }:[0,T]\times \mathbb{R}^d\times \mathcal{P}(U)\times \mathcal{P}(V)\rightarrow \mathcal{S}_{d\times d}$
by 
$$\tilde{\sigma } (t,x, \mu, \nu)\triangleq \left (\int _{U\times V}\sigma (t,x,u,v) \sigma ^T (t,x,u,v)\mu (du) \nu (dv)\right) ^{\frac 12},$$ 
where $()^{\frac 12}$ is the square root of symmetric non-negative definite matrices.
Recall that $\mathcal{P}(U)$ and $\mathcal{P}(V)$ are compact spaces  when endowed with the  (metrizable)  topology of weak convergence. In addition, we have
\begin{Proposition}
The coefficients $\tilde{b}$ and $\tilde{\sigma}$ are jointly continuous, locally  uniformly Lipschitz in $x$ and have  global linear growth (i.e. satisfy items 1-3 from the Standing Assumptions written for  $(\tilde{b},\tilde{\sigma})$). 
\end{Proposition}
Proof: joint continuity is quite simple, since $b$ and $\sigma$ are jointly uniformly continuous on \\
$[0,T]\times \{x\in \mathbb{R}^d||x|\leq K\}\times U \times V$ with modulus of continuity $\omega (K, \cdot)$. Now, if 
$t_n\rightarrow t, $ $x_n\rightarrow x$, $\mu _n\rightarrow \mu $ and $\nu _n\rightarrow n$, we have
$$ \tilde {b}(t_n, x_n, \mu_n, \nu _n)-\tilde{b}(t,x, \mu, \nu)=
 \tilde {b}(t_n, x_n, \mu_n, \nu _n)-\tilde{b}(t,x, \mu_n, \nu_n)+
 \tilde {b}(t, x, \mu_n, \nu _n)-\tilde{b}(t,x, \mu, \nu).$$ Now, as long as $|x_n|\leq K$ we have
 $$|\tilde {b}(t_n, x_n, \mu_n, \nu _n)-\tilde{b}(t,x, \mu_n, \nu_n) |\leq \omega (K, |(t_n,x_n)-(t,x)|)\rightarrow 0.$$
On the other hand, we have $(\mu_n, \nu _n)\rightarrow (\mu,\nu)$ in $\mathcal{P}(U)\times \mathcal{P}(V)$, so
$$ \tilde {b}(t, x, \mu_n, \nu _n)=\int _{U\times V}b(t,x,u,v)\mu_n (du) \nu _n (dv) \rightarrow \int _{U\times V}b(t,x,u,v)\mu (du) \nu (dv)=\tilde{b}(t,x,\mu,\nu)$$
since, for fixed $t,x$ the function $b$ is bounded and continuous  on $U\times V$. Similarly, $\tilde{\sigma}$ is jointly continuous since $\sigma \rightarrow (\sigma)^{\frac 12}$ is a continuous operation over $\mathcal{S}_{d\times d}$.
Now, $|b(t,x,u,v)-b(t,y,u,v)|\leq L(K)|x-y|$, $ \forall \ |x|,|y|\leq K$ and $ \forall \ u,v$ so integrating we get 
 $$|\tilde {b}(t,x,\mu ,\nu)-b(t,y,\mu ,\nu )|\leq L(K) |x-y|\ \forall t, \mu,\nu \ \textrm{and}\ |x|,|y|\leq K.$$ 
 To show that $\tilde{\sigma}$ is locally uniformly Lipschitz in $x$ we have to separate two possible situations. First, the one dimensional case $d=1$ is trivial, and amounts to the simple (consequence of) triangle inequality for $L^2$ norms
\begin{eqnarray*}
\left | \left ( \int _{U\times V}|\sigma (t,x,u,v)|^2\mu (du) \nu (dv)\right )^{\frac 12}-\left ( \int _{U\times V}|\sigma (t,y,u,v)|^2\mu (du) \nu (dv)\right )^{\frac 12} \right| \leq \\
 \left ( \int _{U\times V}|\sigma (t,x,u,v) -\sigma (t,y,u,v)|^2\mu (du) \nu (dv)\right )^{\frac 12} \leq L(K)|x-y|,\ \ |x||y|\leq K.
\end{eqnarray*}
Second, the multidimensional case is actually a well studied problem, of Lipschitz dependence on the square root of a symmetric and non-negative definite matrix on a parameter. For this, we need the extra conditions in the Standing Assumptions item 4. Using Lebesgue's Dominated Convergence Theorem  together with the fact that, as long as $x$ is bounded all quantities involved (including all derivatives in $x$ up to second order) remain bounded, uniformly in $(t,u,v)$, we conclude that, for each fixed $t,\nu, \mu$ the matrix valued mapping
$$\mathbb{R}^d\ni x\rightarrow \int _{U\times V}\sigma (t,x,u,v) \sigma ^T(t,x,u,v)\mu (du) \nu (dv)\in \mathcal{S}_{d\times d}$$
is two times differentiable, and the second derivatives are actually continuous in $x$. In addition we have
$$\frac{\partial ^2 }{\partial x_i x_j}\left (\int _{U\times V}\sigma (t,x,u,v) \sigma ^T(t,x,u,v)\mu (du) \nu (dv) \right )=\int _{U\times V}\left (\frac{\partial^2 }{\partial x_ix_j} \sigma (t,x,u,v) \sigma ^T(t,x,u,v)\right)\mu (du) \nu (dv).
$$
From item 4 in the Standing Assumption, we conclude that 
$$\left |\frac{\partial ^2 }{\partial x_i x_j}\left (\int _{U\times V}\sigma (t,x,u,v) \sigma ^T(t,x,u,v)\mu (du) \nu (dv) \right )\right|\leq C'(K),$$
$\forall t,u,v$ and $|x|\leq K$, for some $C'(K)<\infty$. Using Theorem 5.2.3  page 132 in \cite{MR2190038}, we conclude that, for fixed $t, \mu,\nu$ the mapping 
$$\{x|\ |x|\leq K\} \ni x\rightarrow  \left (\int _{U\times V}\sigma (t,x,u,v) \sigma ^T (t,x,u,v)\mu (du) \nu (dv)\right) ^{\frac 12}\in \mathcal{S}_{d\times d}$$
is Lipschitz  and the Lipschitz constant does not depend on $t,\mu, \nu$. In other words, $\tilde{\sigma}$ is uniformly locally Lipschitz in $x\in \mathbb{R}^d$.   The global linear growth condition for both $\tilde{b}$ and $\tilde{\sigma}$ is obvious.     $\diamond$

  Now, the whole proof of Theorem \ref{main} hinges on the very simple observation, that elementary mixed feedback strategies $\mu \in \mathcal{U}^m(s)$ and $\nu \in \mathcal{V}^m(s)$ actually become elementary  pure  strategies for the auxiliary game
with state equation

\begin{equation}\label{state-aux}
\left \{
\begin{array}{ll}
d\tilde{X}_t=\tilde{b}(t,\tilde{X}_t,  \mu , \nu )dt+\tilde{\sigma} (t, \tilde{X}_t,\mu , \nu )d \tilde{W}_t,\ s \leq t\leq T\\
\tilde{X}_{s}=x\in \mathbb{R}^d,
\end{array}
\right .
\end{equation} 
where $\tilde{W}$ is $d$-dimensional (same as the state) Brownian Motion. In other words, using tilde in an obvious way, we have
$$\tilde{\mathcal{U}}^p(s)=\mathcal{U}^m(s), \tilde{\mathcal{V}}^p(s)=\mathcal{V}^m(s).$$
 We apply Proposition \ref{prop:state-eq-simple} (replacing $U$ by $\mathcal{P}(U)$ does not make any difference) to the auxiliary equation \eqref{state-aux},  for a \emph{fixed}   Brownian motion $(\tilde{W}_t)_{0\leq t\leq T}$  on $(\tilde{\Omega}, \tilde{\mathcal{F}}, \tilde{\mathbb{P}})$ with natural
 filtrations $$\tilde{\mathcal{F}}^s_t=\sigma (\tilde{W}_u-\tilde{W}_s, s\leq u\leq t)\vee \mathcal{N}(\tilde{\mathbb{P}}, \tilde{ \mathcal{F}})\ \  \textrm{for}\ \ s\leq t\leq T.$$ 
We obtain that  
\begin{enumerate}
\item equation \eqref{state-aux} has a unique strong solution $\tilde{X}^{s,x,\mu,\nu}$
\item equation \eqref{state-aux} satisfies the path-wise uniqueness property.
\end{enumerate}
Again, using the celebrated Yamada-Watanabe result in \cite{KS88}, existence and  path-wise uniqueness implies existence and uniqueness in law of a weak solution to \eqref{state-aux}. Weak existence and uniqueness in law is equivalent,  following  \cite{MR2190038},  to the existence and uniqueness of the local Martingale Problem associated with \eqref{state-aux}. However, the local Martingale Problem associated to \eqref{state-aux}  is \emph{identical} to the  the local martingale problem in Definition \ref{mart}. Combining, we have the following lemma, that contains actually more than the proof of Proposition \ref{prop:well-posed}.
\begin{Lemma} Fix $s,x$ and $\mu \in \mathcal{U}^m(s)$, $\nu \in \mathcal{V}^m(s)$.  The law of $\tilde{X}^{s,x,\mu, \nu}$ i.e.
$\mathbb{Q}^{s,x, \mu ,\nu}\triangleq \tilde{\mathbb{P}}\circ \left (\tilde{X}^{s,x, \mu, \nu}_{\cdot}\right )^{-1}$,
is the unique solution of the martingale problem in Definition \ref{mart} (associated with the original state equation with initial conditions $s,x$ and elementary continuous time mixed strategies $\mu, \nu$).
\end{Lemma}
From the above Lemma, we can also conclude easily that
$$J(s,x,\mu, \nu)=\tilde{\mathbb{E}}[g(\tilde{X}^{s,x, \mu, \nu}_T)].$$
Therefore, the values of the original game  satisfy:
$$
V^-(s,x)= \sup _{\mu  \in \mathcal{U}^m(s)} \inf _{\nu \in \mathcal{V}^m(s)} \tilde{\mathbb{E}}[g(\tilde{X}^{s,x, \mu, \nu}_T)]\leq 
\inf _{\nu \in \mathcal{V}^m(s)}  \sup _{\mu  \in \mathcal{U}^m(s)} \tilde{\mathbb{E}}[g(\tilde{X}^{s,x, \mu, \nu}_T)]=V^+(s,x),$$
which is all the same as 
$$
V^-(s,x)= \sup _{\mu  \in \tilde{\mathcal{U}}^p(s)} \inf _{\nu \in \tilde{\mathcal{V}}^p(s)} \tilde{\mathbb{E}}[g(\tilde{X}^{s,x, \mu, \nu}_T)]\leq 
\inf _{\nu \in \tilde{\mathcal{V}}^p(s)}  \sup _{\mu  \in \tilde{\mathcal{U}}^p(s)} \tilde{\mathbb{E}}[g(\tilde{X}^{s,x, \mu, \nu}_T)]=V^+(s,x).$$
In words, we have translated the original game over mixed strategies in a game  over pure strategies  in strong formulation, for the new state $\tilde{X}$. The new game, however, satisfies the Isaacs condition, and the (common) Isaacs equation coincides with the randomized Isaacs equation \eqref{eq:Isaacs} for the original game. The coefficients $\tilde{b}$ and $\tilde{\sigma}$ satisfy the standing assumptions in \cite{sirbu}. Therefore, we can use Theorem 4.1 in \cite{sirbu} (which does apply if one replaces $U$ by $\mathcal{P}(U)$ and $V$ by $\mathcal{P}(V)$) to conclude that $V^-=V^+$ and the common value is the unique bounded continuous viscosity solution of 
\eqref{eq:Isaacs}.  In addition, for each $\rho \in \mathbb{B}^s$ we have the (DPP)
$$V(s,x)= \sup _{\mu \in \mathcal{U}^m(s)}\inf _{\nu \in \mathcal{V}^m(s)}  \tilde{\mathbb{E}}\left [V(\rho (\tilde{X}^{s,x,\mu, \nu}_{\cdot}),\tilde{ X}^{s,x,\mu, \nu}_{ \rho (\tilde{X}^{s,x,u,v}_{\cdot}   )} \right]
=\inf _{\nu \in \mathcal{V}^m(s)}   \sup _{\mu \in \mathcal{U}^m(s)} \tilde{\mathbb{E}}\left [V(\rho (\tilde{X}^{s,x,\mu, \nu}_{\cdot}),\tilde{ X}^{s,x,\mu, \nu}_{ \rho (\tilde{X}^{s,x,u,v}_{\cdot}   )} \right].$$
Recalling that $\mathbb{Q}^{s,x, \mu ,\nu}\triangleq \tilde{\mathbb{P}}\circ \left (\tilde{X}^{s,x, \mu, \nu}_{\cdot}\right )^{-1},$ we obtain the DPP for the original problem, i.e.
$$V(s,x)= \sup _{\mu \in \mathcal{U}^m(s)} \inf _{\nu \in \mathcal{V}^m(s)}  \mathbb{E}^{\mathbb{Q}^{s,x,\mu,\nu }}[V (\rho(y), y(\rho (y)))]=\inf _{\nu \in \mathcal{V}^m(s)}  \sup _{\mu \in \mathcal{U}^m(s)}  \mathbb{E}^{\mathbb{Q}^{s,x,\mu ,\nu }}[V (\rho(y), y(\rho (y)))], \forall \rho \in \mathbb{B}^s.$$
\section{Conclusions}
We define, symmetrically, a continuous-time game over   continuous-time mixed feedback strategies.  These are feedback counterparts to relaxed/generalized open-loop controls in the sense of \cite{young}.
The state equation is well posed in the sense of the  local  Martingale Problem of Stroock and Varadhan  as long as we restrict the game to the class of elementary strategies.  We show that the game formulated over such martingale solutions  has a value, which is the unique continuous viscosity solution of the randomized Isaacs equation. To the best of our knowledge, such result is not present in the literature. 

 In our view, the present result does not come in competition with the   very recent work \cite{buc-li-q}, but it is actually complementary.  Beyond 
the different interpretation of  the the state equation,  \cite{buc-li-q}
 show that the lower and upper values  ($W^{\pi}$ and $U^{\pi}$) of the game discretized over  a fixed time grid $\pi$ converge to the same value as $|\pi|\rightarrow 0$. In other words, they first discretize the game, and then consider the limit.
We consider directly limiting strategies, by allowing continuous-time mixing (that would be relaxed controls if they were open-loop). We can show that, with this limiting definition of strategies, the original continuous-time game has a value. The intuition tells us clearly that one of the two routes has to be considered to obtain a value. If the weaker player (the one outside the inf sup/ sup inf) is restricted to randomize over a discrete-time grid (consisting of deterministic times or stopping-rules), the changes of actions of the stronger player should be restricted to a similar (if not the same) time grid. If the stronger player can change actions much faster than the weaker player can randomize, the stronger player would learn immediately the result of the randomization of the weaker player, and therefore take advantage of it, as if the weaker player's strategy were pure. In other words, one would still obtain the (different) upper and lower values of the game over pure strategies.

\bibliographystyle{amsalpha} 
\providecommand{\bysame}{\leavevmode\hbox to3em{\hrulefill}\thinspace}
\providecommand{\MR}{\relax\ifhmode\unskip\space\fi MR }
\providecommand{\MRhref}[2]{%
  \href{http://www.ams.org/mathscinet-getitem?mr=#1}{#2}
}
\providecommand{\href}[2]{#2}

\end{document}